\documentclass{elsarticle}

\usepackage{amsmath, amssymb, cases}

\usepackage{bm}

\newcommand{\secref}[1]{\S#1}

\newcommand{\Vector}[1]{\bm{#1}}

\newcommand{\Triangle}{\triangle}
\newcommand{\Cantor}{\mathcal{C}}

\newcommand{\Amount}[1]{|#1|}

\newcommand{\R}{\mathbb{R}} % Real numbers
\newcommand{\N}{\mathbb{N}} % Natural numbers

\DeclareMathOperator{\Span}{span}

\begin{document}
\title{Fractal patterns related to dividing coins}
\author{Ken Yamamoto}
\ead{yamamoto@phys.chuo-u.ac.jp}
\address{Department of Physics, Faculty of Science and Engineering, Chuo University, Kasuga, Bunkyo, Tokyo 112-8851, Japan}

\begin{abstract}
The present paper formulates and solves a problem of dividing coins.
The basic form of the problem seeks
the set of the possible ways of dividing
coins of face values $1,2,4,8,\ldots$ between three people.
We show that this set possesses a nested structure like the Sierpinski-gasket fractal.
For a set of coins with face values power of $r$,
the number of layers of the gasket becomes $r$.
A higher-dimensional Sierpinski gasket is obtained if the number of people is more than three.
In addition to Sierpinski-type fractals,
the Cantor set is also obtained in dividing an incomplete coin set
between two people.
\end{abstract}


\maketitle

\section{Introduction}\label{sec1}
Money is a mathematical system so familiar to us.
It is a very good instance of combinatorics and discrete mathematics.

Telser \cite{Telser} considered the problem of optimal currency,
and deduced that
the optimal currency system consists of denominations
with face values of $1, 3, 9,\ldots, 3^{m-1}$.
This result comes from the problem of Bachet
which seeks the smallest number of weights
so that they can weigh any integer quantity on a two-pan balance.
Many actual currency systems take their average multiples close to three:
$2.8$ by Wynne \cite{Wynne} and $2.60$ by Tschoegl \cite{Tschoegl}.
This fact is surprising
because most currencies are established based on the decimal system,
which is not compatible with Telser's power-of-three theory,
and because individual customs and culture of each country
are reflected to its currency.

One of the classical and famous problems is the Frobenius coin problem,
which seeks the largest amount of money
that cannot be made using only coins of specified denominations
\cite{RamirezAlfonsin2005}.
The solution to this problem is called the Frobenius number.
For only two types of coins of denominations $a_1$ and $a_2$
which are relatively prime,
the Frobenius number is given by $a_1a_2-a_1-a_2$ \cite{Sylvester}.
However, closed expressions are not known for more than two types of coins,
and the problem was found to be NP-hard \cite{RamirezAlfonsin1996}.

Another classical problem about a coin is
the change-making problem.
It asks
how a given amount of money can be made with the least number of coins
of given denominations.
This problem is a variant of the knapsack problem \cite{Kellerer},
and is also an NP-hard problem \cite{Lueker}.
The greedy algorithm \cite{Magazine} and other heuristic methods
such as dynamic programming \cite{Martello}
give the optimal solution in some cases.
A coin system is said to be {\it canonical}
if the greedy algorithm works correctly \cite{Pearson},
and almost all currencies in the world are arranged to be canonical.
Even with a canonical coin system, which is easy to study,
the change-making process possesses a rich mathematical structure.
If we repeatedly pay money
so that the number of coins in the purse after each payment is minimized,
a fractal pattern is obtained
from a sequence of change amounts \cite{Yamamoto2012, Yamamoto2013}.

The present article develops a generation mechanism of a fractal set
in the form of dividing coins between people.
The problem is very simple in appearance,
but the possible ways of division form a nontrivial fractal pattern
like the Sierpinski gasket and Cantor set.
The basic situation, discussed in \secref{\ref{sec2}}, is that three people divide coins
of face values $1,2,4,\ldots,2^{m-1}$,
and the Sierpinski gasket appears as the attractor.
The emergence of the Sierpinski gasket has been previously reported \cite{Yamamoto},
but in the present article we rigorously formulate 
the fractal structure and diverse generalizations
in the framework of iterated function systems.
In addition to the Sierpinski gasket, 
the Cantor set also appears in a specific case.

The notion of fractal was invented by Mandelbrot \cite{Mandelbrot}
to measure the morphology of natural objects.
Fractal theory has been applied to various phenomena
such as in physics, economics, and biology.
In naive description,
a fractal object is created by defining a starting shape (an {\it initiator})
and replacing each part
with another shape called a {\it generator}, ad infinitum \cite{Feder}.
Mathematically, a fractal is defined
using an iterated function system (IFS for short).
Here we briefly review construction of the Sierpinski gasket and Cantor set
(see Barnsley \cite{Barnsley} for detail).
Let $\Vector{p}_1,\Vector{p}_2,\Vector{p}_3$ be three points
in $\R^d$ which are not collinear.
Let three contraction maps $f_1,f_2,f_3:\R^d\to\R^d$ be given by
$f_i(\Vector{x})=(\Vector{x}+\Vector{p}_i)/2$ for $i=1,2,3$.
The set of contraction maps $\{f_1,f_2,f_3\}$ is called an IFS.
We let $\mathcal{H}(\R^d)$ 
denote the collection of all nonempty compact subsets of $\R^d$.

It is well known that $\mathcal{H}(\R^d)$ is a complete metric space with
Hausdorff metric \cite{Barnsley}.
The IFS $\{f_1,f_2,f_3\}$ defines a map
$F:\mathcal{H}(\R^d)\ni K
	\mapsto f_1(K)\cup f_2(K)\cup f_3(K)\in\mathcal{H}(\R^d)$
which is a contraction on $\mathcal{H}(\R^d)$.
The Sierpinski gasket $\Triangle(\in\mathcal{H}(\R^d))$
is defined as the unique fixed point of $F$ which obeys
$\Triangle=F(\Triangle)=f_1(\Triangle)\cup f_2(\Triangle)\cup f_3(\Triangle)$;
existence and uniqueness are guaranteed
by the contraction mapping theorem (or the Banach fixed-point theorem).
The fixed point $\Triangle$ also satisfies
$\Triangle=\lim_{m\to\infty}F^m(K)$ for any $K\in\mathcal{H}(\R^d)$;
in this sense, $\Triangle$ is called the {\it attractor} of the IFS.
The three points $\Vector{p}_1$, $\Vector{p}_2$, and $\Vector{p}_3$
are at the three corners of $\Triangle$.
In a similar way,
the Cantor set between $\Vector{p}_1$ and $\Vector{p}_2$ is the attractor of an IFS $\{g_1,g_2:\R^d\to\R^d\}$
where $g_i(\Vector{x})=(\Vector{x}+\Vector{p}_i)/3$ for $i=1,2$.

\section{Formulation of the problem and a basic solution} \label{sec2}
This paper treats of a problem related to
dividing a set of coins by some ``players''.
In order to specify the set of coins,
we introduce a notation
\[
\begin{pmatrix}
v_1 & \cdots & v_m\\ c_1 & \cdots & c_m
\end{pmatrix},
\]
where $v_1,\ldots,v_m$ is the list of face values of the coins,
and $c_i$ is the number of coins of face value $v_i$.
For a coin set
$S=\begin{pmatrix}v_1 & \cdots & v_m\\ c_1 & \cdots & c_m\end{pmatrix}$,
we define $\Amount{S}$ as the total amount of coins in $S$, given by
\[
\Amount{S}:=\sum_{i=1}^m v_i c_i.
\]

For example, the set of coins
\[
S=\begin{pmatrix}1 & 10\\ 2 & 1\end{pmatrix}
\]
consists of two pennies (1 cent coins) and one dime (10 cent coin).
When two players $A$ and $B$ divide these coins,
the division of these coins is expressed 
by a pair $(n_A, n_B)$ of money amounts that $A$ and $B$ respectively receive.
There are six ways to divide the coins in this case:
\[
(n_A,n_B)=(0,12), (1,11), (2,10), (10,2), (11,1), (12,0).
\]
Note that we admit cases in which some players receive no coins.
In this paper, we mainly focus on the coin-dividing problem between three players.

We start from the most basic form of the problem.
Let us suppose a ``binary'' currency system
\[
S_{2,m}=\begin{pmatrix}1 & 2 & 4 &\cdots& 2^{m-1}\\
1& 1& 1&\cdots& 1\end{pmatrix},
\]
and three players $A$, $B$, and $C$ divide these coins.
Each way of division is represented
by lining each player's share as the triplet $(n_A,n_B,n_C)$,
which defines a point in a three-dimensional space.
We study the structure of the possible points $(n_A,n_B,n_C)$.

This problem can be described by weights.
If there are $m$ weights of $1, 2,\ldots, 2^{m-1}$ grams,
and one separates them into three groups (admitting one or two null groups),
then we study the shape of the possible groupings.
The equivalence of coins and weights resembles
Telser's theory of optimal currency
(see the beginning of \secref{\ref{sec1}}).

\begin{figure}[!tb]\centering
\mbox{\raisebox{40mm}{(a)}}
\mbox{\raisebox{10mm}{\includegraphics[clip]{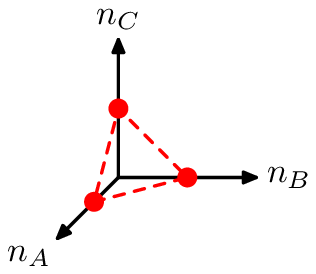}}}
\hspace{5mm}
\mbox{\raisebox{40mm}{(b)}}
\includegraphics[clip]{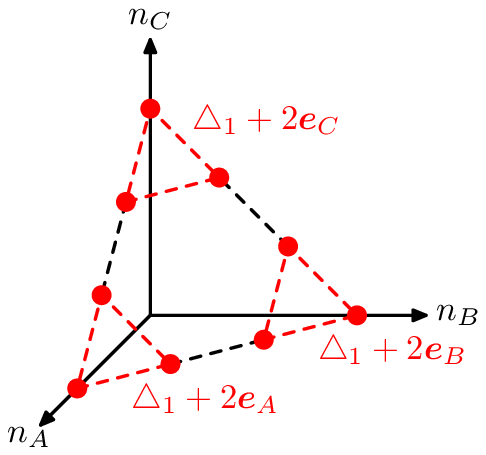}
\hspace{5mm}
\mbox{\raisebox{40mm}{(c)}}
\includegraphics[clip]{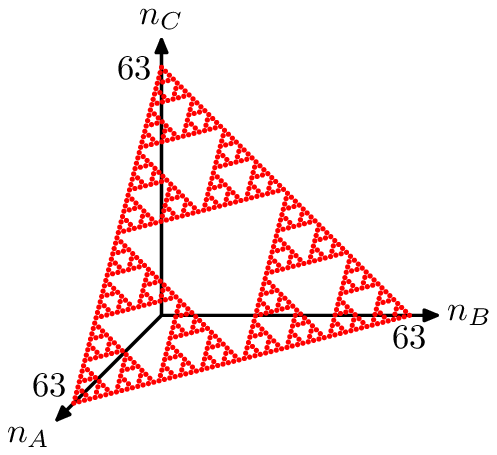}
\caption{
(a) The set $\Triangle_1$ consists of three points
corresponding to the basis vectors
$\Vector{e}_A$, $\Vector{e}_B$ and $\Vector{e}_C$,
which form an equilateral triangle (the dashed lines).
(b) $\Triangle_2$ is composed of three copies of $\Triangle_1$
placed at the three corners.
(c) $\Triangle_7$ looks like the Sierpinski gasket.
}
\label{fig1}
\end{figure}

Formally, an orthonormal basis $\{\Vector{e}_A,\Vector{e}_B,\Vector{e}_C\}$
is taken
so that $(n_A, n_B, n_C)=n_A\Vector{e}_A+n_B\Vector{e}_B+n_C\Vector{e}_C$.
By definition each division satisfies $n_A+n_B+n_C=2^m-1$
and $n_A,n_B,n_C\ge0$.
Hence, all the possible points $(n_A,n_B,n_C)$ are confined
in a triangular area spanned by
$(2^m-1)\Vector{e}_A$, $(2^m-1)\Vector{e}_B$, and $(2^m-1)\Vector{e}_C$.
Let $\Triangle_m (\subset\R^3)$ denote the set of the points
of possible division $(n_A,n_B,n_C)$.
For the simplest case $m=1$,
three players divide only one coin of face value 1.
The possible ways of division are
$(n_A,n_B,n_C)=(1,0,0)$, $(0,1,0)$, and $(0,0,1)$,
and they form an equilateral triangle in the three dimensional space
(Fig.~\ref{fig1}(a)):
\[
\Triangle_1=\{(1,0,0), (0,1,0), (0,0,1)\}
=\{\Vector{e}_A, \Vector{e}_B, \Vector{e}_C\}.
\]

Next we take into account the coin of face value 2 to advance to $m=2$.
If the player $A$ receives this coin,
the set of possible division is obtained
by considering who receives the coin of face value 1.
This set is written as
$\{(3,0,0), (2,1,0),(2,0,1)\}=\Triangle_1+2\Vector{e}_A$,
which is translation of $\Triangle_1$.
Similar expressions hold for $B$ and $C$,
so $\Triangle_2$ is given by
\[
\Triangle_2=(\Triangle_1+2\Vector{e}_A)\cup
		(\Triangle_1+2\Vector{e}_B)\cup
		(\Triangle_1+2\Vector{e}_C).
\]
That is, $\Triangle_2$ consists of three copies of $\Triangle_1$,
each of which is placed at one of the three corners
(see Fig.~\ref{fig1}(b) for reference).
Similarly, by considering who receives the coin of face value $2^{m-1}$,
$\Triangle_m$ is inductively given by
\begin{equation}
\Triangle_m=(\Triangle_{m-1}+2^{m-1}\Vector{e}_A)\cup
		(\Triangle_{m-1}+2^{m-1}\Vector{e}_B)\cup
		(\Triangle_{m-1}+2^{m-1}\Vector{e}_C).
\label{eq:induction}
\end{equation}
For each $m\ge2$,
three copies of $\Triangle_{m-1}$ are placed triangularly
to form $\Triangle_m$.
By setting $\Triangle_0=\{(0,0,0)\}$ for convenience,
Eq.~\eqref{eq:induction} holds also for $m=1$.

Intuitively,
$\Triangle_m$ achieves self-similarity as $m$ increases.
Actually, $\Triangle_7$, shown in Fig.~\ref{fig1} (c),
is almost indistinguishable from the genuine Sierpinski gasket.
The connection between $\Triangle_m$ and the Sierpinski gasket
is formulated as follows.
We set three contraction maps $f_A,f_B,f_C:\R^3\to\R^3$ given by
\begin{equation}
f_P(\Vector{x})=\frac{\Vector{x}+\Vector{e}_P}{2}\quad (P=A,B,C).
\label{eq:fP}
\end{equation}
As stated in \secref{\ref{sec1}}, the attractor $\Triangle\in\mathcal{H}(\R^3)$
of an IFS $\{f_A,f_B,f_C\}$ is the Sierpinski gasket
spanned between $\Vector{e}_A$, $\Vector{e}_B$, and $\Vector{e}_C$.
On the other hand,
multiply $2^{-m}$ to Eq.~\eqref{eq:induction} to get
\begin{align*}
2^{-m}\Triangle_m
&=\left(\frac{2^{-(m-1)}\Triangle_{m-1}+\Vector{e}_A}{2}\right)\cup
  \left(\frac{2^{-(m-1)}\Triangle_{m-1}+\Vector{e}_B}{2}\right)\cup
  \left(\frac{2^{-(m-1)}\Triangle_{m-1}+\Vector{e}_C}{2}\right)\\
&=f_A(2^{-(m-1)}\Triangle_{m-1})\cup
	f_B(2^{-(m-1)}\Triangle_{m-1})\cup f_C(2^{-(m-1)}\Triangle_{m-1})\\
&=F(2^{-(m-1)}\Triangle_{m-1}),
\end{align*}
where $F$ is the contraction map on $\mathcal{H}(\R^3)$ induced from
the IFS.
That is, adding the $(m+1)$st coin to the coin set $S_{2,m}$ 
is directly represented as the action of $F$.
Hence $2^{-m}\Triangle_m=F^m(\Triangle_0)$.
Because $\Triangle_0=\{(0,0,0)\}\in\mathcal{H}(\R^3)$,
it follows from the contraction mapping theorem that
\[
\Triangle=\lim_{m\to\infty}F^m(\Triangle_0)=\lim_{m\to\infty}2^{-m}\Triangle_m.
\]
This is the result
which connects the coin-dividing problem and the Sierpinski gasket.

\section{Generalization of the coin system}
In this section we consider a generalized coin-dividing problem where
the coin system is not powers of two.
We give the set of coins
\[
S_{r,m}=\begin{pmatrix}1 & r & r^2 &\cdots& r^{m-1}\\
	r-1& r-1& r-1&\cdots& r-1\end{pmatrix},
\]
where $r$ is an integer equal to or greater than two,
and $\Amount{S_{r,m}}=r^m-1$.

Before working on the problem of $S_{r,m}$,
we revise Eq.~\eqref{eq:induction} to derive a useful expression.
By means of the result
$\Triangle_1=\{\Vector{e}_A, \Vector{e}_B, \Vector{e}_C\}$,
Eq.~\eqref{eq:induction} is rewritten as
\[
\Triangle_m
=\bigcup_{\Vector{q}\in\Triangle_1}(\Triangle_{m-1}+2^{m-1}\Vector{q})
\]
for $m\ge2$.
A similar formula is valid for $S_{r,m}$, and hence
the set $\Triangle_{r,m}$ of the possible points corresponding to $S_{r,m}$
is inductively given by
\begin{align}
\Triangle_{r,1}
&=\{q_A\Vector{e}_A+q_B\Vector{e}_B+q_C\Vector{e}_C|
	q_A,q_B,q_C\in\N\cup\{0\}, q_A+q_B+q_C =r-1\}, \label{eq:triangle_r1}\\
\Triangle_{r,m}&=\bigcup_{\Vector{q}\in\Triangle_{r,1}}
	(\Triangle_{r,m-1}+r^{m-1}\Vector{q}).
\label{eq:triangle_rm}
\end{align}
The set $\Triangle_{r,1}$ consists of $r(r+1)/2$ points
which form a triangle,
and similarly
$\Triangle_{r,m}$ consists of $r(r+1)/2$ copies of $\Triangle_{r,m-1}$
arranged triangularly.

\begin{figure}[!t]\centering
\mbox{\raisebox{40mm}{(a)}}
\includegraphics[height=40mm]{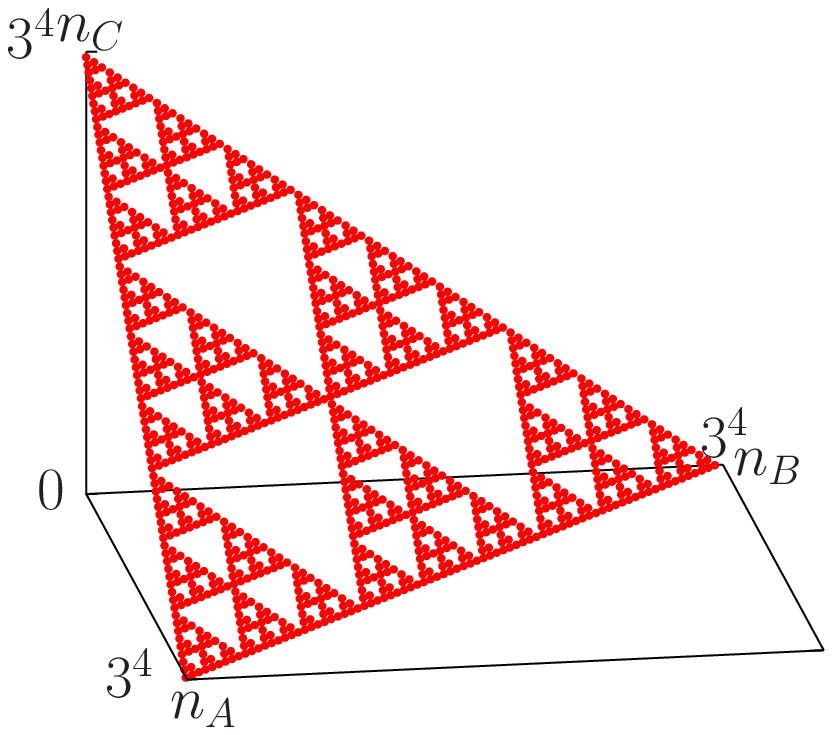}
\mbox{\raisebox{40mm}{(b)}}
\includegraphics[height=40mm]{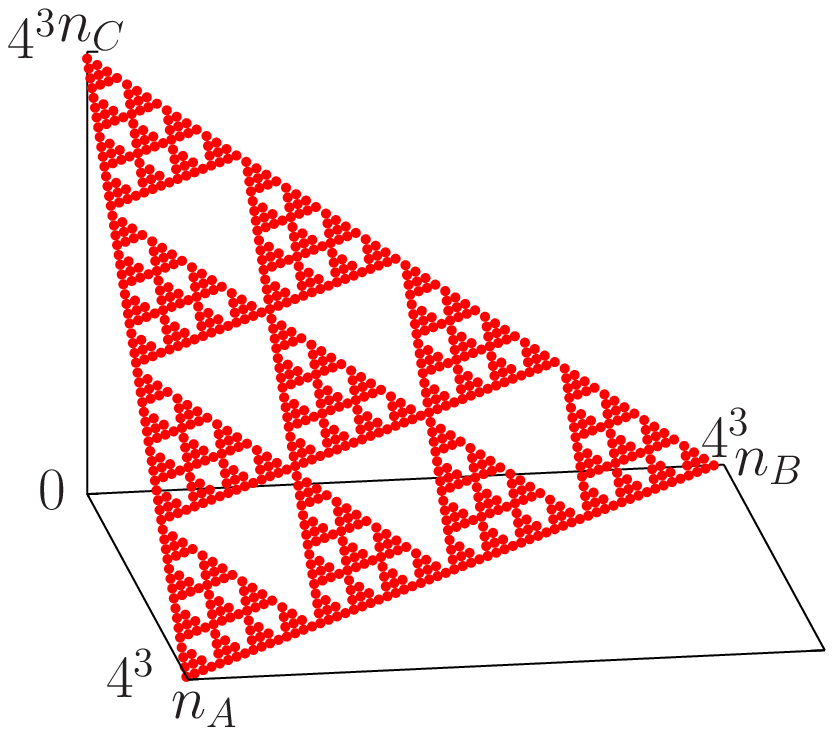}
\mbox{\raisebox{40mm}{(c)}}
\includegraphics[height=40mm]{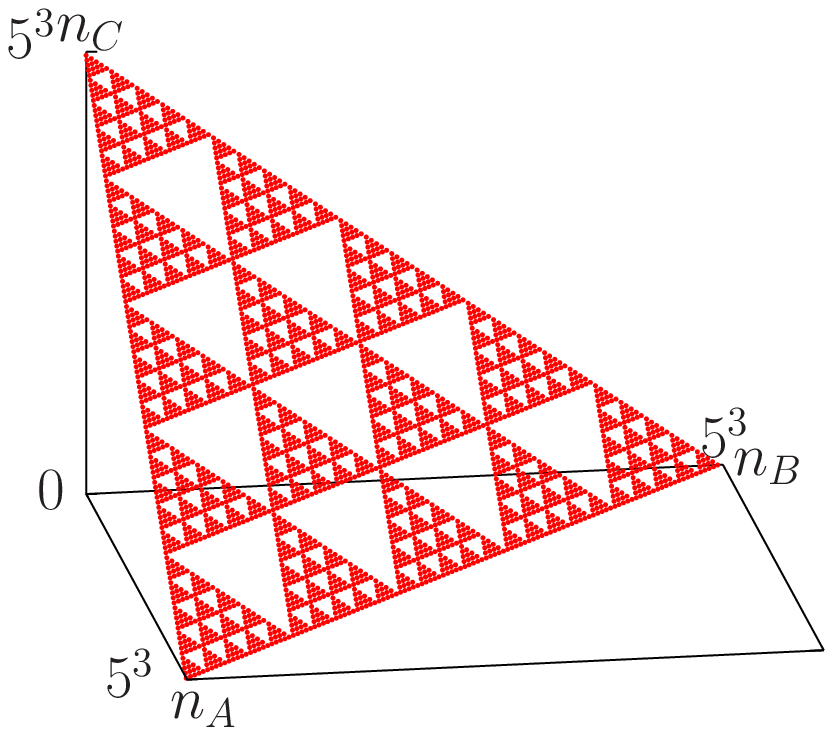}
\caption{
The sets of possible points $(n_A,n_B,n_C)$ of the coin set
$S_{3,4}$ (a), $S_{4,3}$ (b), and $S_{5,3}$ (c).
They resemble the Sierpinski gasket,
but the number of layers is equal to $r$.
}
\label{fig2}
\end{figure}

By multiplying $r^{-m}$ to Eq.~\eqref{eq:triangle_rm},
we have
\[
r^{-m}\Triangle_{r,m}=\bigcup_{\Vector{q}\in\Triangle_{r,1}}
	f_{\Vector{q}}(r^{-(m-1)}\Triangle_{r,m-1}),
\]
where $f_{\Vector{q}}:\R^3\ni\Vector{x}\mapsto(\Vector{x}+\Vector{q})/r\in\R^3$
is a contraction map.
By the contraction mapping theorem,
$r^{-m}\Triangle_{r,m}$ converges to an attractor
of IFS $\{f_{\Vector{q}}\vert \Vector{q}\in\Triangle_{r,1}\}$, as $m\to\infty$.
This attractor is the Sierpinski gasket with $r$ layers.
We depict the coin-dividing plot corresponding to
$r=3$, $4$, and $5$ in Fig.~\ref{fig2}.
The plot is essentially the same as a class of Pascal-Sierpinski gaskets
if $r$ is prime---
the Pascal-Sierpinski gasket of order $p$ is obtained
by coloring the number not divisible by $p$ in Pascal's triangle \cite{Holter}.

More generally, we further generalize the coin set
\[
S_{r,c,m}=\begin{pmatrix} 1 & r & \cdots & r^{m-1}\\
	c & c & \cdots & c\end{pmatrix},
\]
and the corresponding set $\Triangle_{r,c,m}$ of possible division of coins.
Obviously, the point set $\Triangle_{r,c,m}$ lies on the plane
given by $n_A+n_B+n_C=\Amount{S_{r,c,m}}=c(r^m-1)/(r-1)$.
Furthermore, the coin-dividing problem is symmetric with respect to
the three players $A$, $B$, and $C$.
Hence $\Triangle_{r,c,m}$ is invariant
under the action of any permutation of $(A,B,C)$.
That is, $\Triangle_{r,c,m}$ has left-right symmetry
and three-fold rotational symmetry.

\begin{figure}[!t]\centering
\mbox{\raisebox{40mm}{(a)}}
\includegraphics[height=40mm]{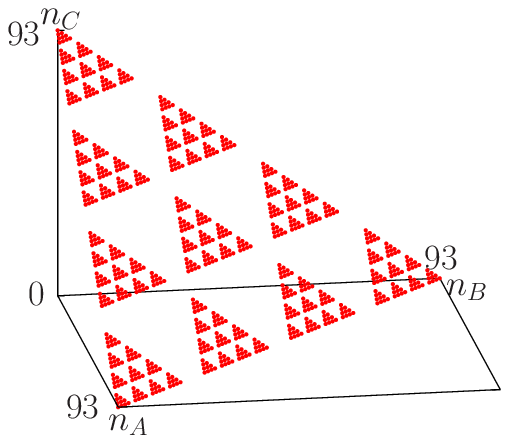}
\mbox{\raisebox{40mm}{(b)}}
\includegraphics[height=40mm]{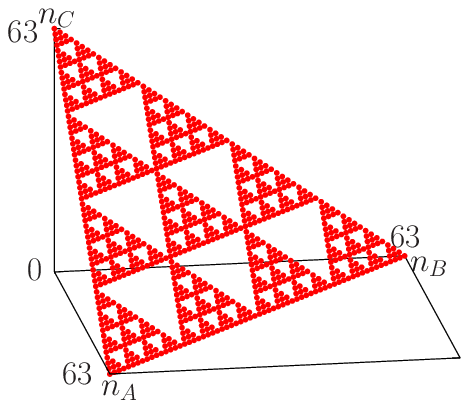}
\mbox{\raisebox{40mm}{(c)}}
\includegraphics[height=40mm]{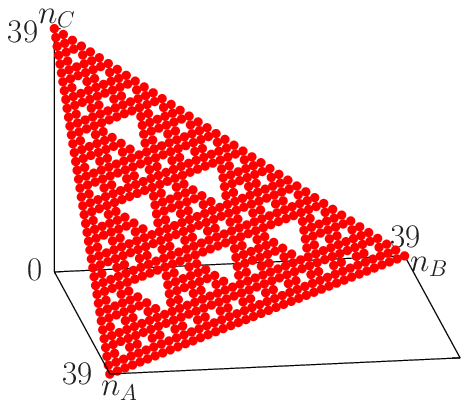}
\caption{
Illustration of $\Triangle_{r,c,m}$ with $c=3$, $m=3$, and $r=5,4,3$.
Three types of fractals appear according to $c$ and $r-1$.
(a) $r=5$ ($c<r-1$). The pattern is totally disconnected.
(b) $r=4$ ($c=r-1$). The pattern is finitely ramified.
(c) $r=3$ ($c>r-1$). The pattern is infinitely ramified.
}
\label{fig3}
\end{figure}

The inductive formula of $S_{r,c,m}$ is
\[
\Triangle_{r,c,m}
=\bigcup_{\Vector{q}\in\Triangle_{c+1,1}}
	(\Triangle_{r,c,m-1}+r^{m-1}\Vector{q}),
\]
where $\Triangle_{c+1,1}$ is given by Eq.~\eqref{eq:triangle_r1}.
The structure of the attractor is classified into three cases
according to the values of $r$ and $c$.
\begin{enumerate}
\item If $c<r-1$, some money amounts less than $\Amount{S_{r,c,m}}$ 
cannot be made by picking out coins from $S_{r,c,m}$.
In fact, the money amount $r-1$ needs $r-1$ coins of face value 1,
but there are only $c$ ($<r-1$) coins of face value 1 in $S_{r,c,m}$.
Reflecting this, the attractor becomes totally disconnected.
\item If $c=r-1$, the attractor is the Sierpinski gasket with $r$ layers.
It belongs to the class of finitely-ramified fractals,
which means that any subset of the fractal can be disconnected
by removing a finite number of points.
\item If $c>r-1$, some money amounts can be made by more than one way;
for example, the money amount $r$ can be made by
one coin of face value $r$ or by $r$ coins of face value 1.
The attractor becomes infinitely ramified.
\end{enumerate}
Figure~\ref{fig3} shows the structural difference
by the magnitude of $c$ and $r-1$.

The Hausdorff dimension of the corresponding attractor 
is easily derived for $c\le r-1$.
The whole pattern can be decomposed
into the $(c+1)(c+2)/2$ smaller copies of itself scaled by a factor $r$,
so the similarity dimension $D_{r,c}$, depending on $r$ and $c$, is given by
\[
D_{r,c}=\frac{\ln((c+1)(c+2)/2)}{\ln r}.
\]
Since $(c+1)(c+2)/2<r^2$ for any $r\ge2$ and $c\le r-1$,
we have $D_{r,c}<2$.
In particular, $D_{2,1}=\ln3/\ln2$ ($\approx1.585$)
is the Hausdorff dimension of the ordinary Sierpinski gasket.
We recall here that
$(c+1)(c+2)/2$ is the number of contraction maps of the IFS,
and $1/r$ is the contraction ratio of each map in the IFS.
For $c>r-1$, the dimension is not calculated easily because of overlap.
The Hausdorff dimension of an overlapping Sierpinski gasket
has been derived only in restricted situations \cite{Broomhead}.

\begin{figure}[t!]\centering
\includegraphics[clip]{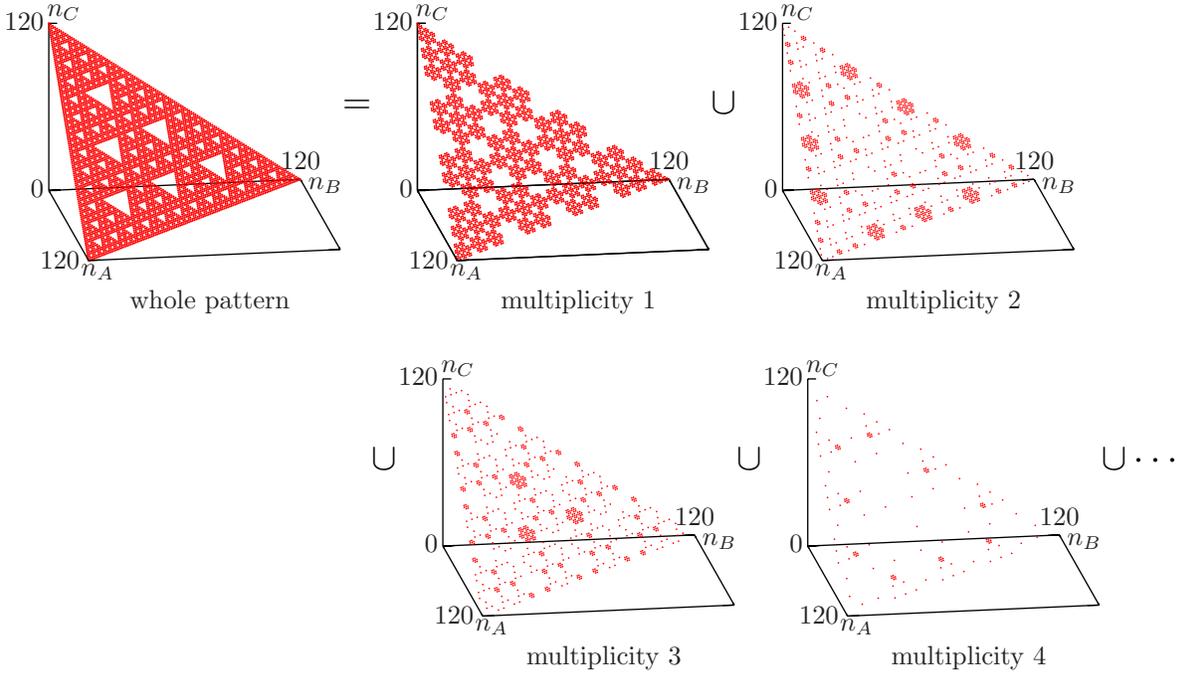}
\caption{
Decomposition of $\Triangle_{3,3,4}$ according to the multiplicity of each point.
The subsets of multiplicity 1, 2, 3, and 4 are shown.
}
\label{fig4}
\end{figure}

In the case of $c>r-1$, a point of $\Triangle_{r,c,m}$ can correspond to more than one ways of coin division.
We can decompose $\Triangle_{r,c,m}$ into subsets according to this multiplicity.
Figure \ref{fig4} shows an example of this decomposition where $(c,r,m)=(3,3,4)$.
The highest multiplicity in this case is nine,
so $\Triangle_{r,c,m}$ splits into nine subsets,
each of which has a fractal shape.
A subset becomes sparse in high multiplicity,
and we only present the lowest four subsets in the figure.
We expect that this multi-level structure gives insight to analyses of an overlapping Sierpinski gasket.

\begin{figure}[t!]\centering
\includegraphics[clip]{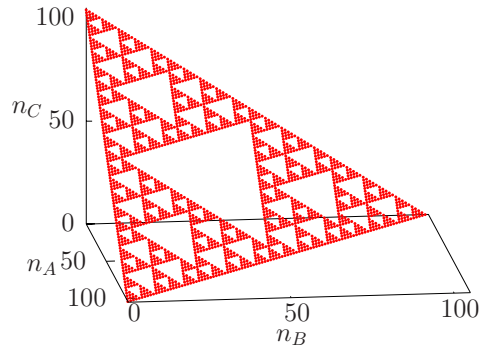}
\caption{
The coin-dividing set of the coin $S_{\mathrm{cent}}$.
A hierarchical structure breaks a little.
}
\label{fig5}
\end{figure}

For reference,
we consider the coin system of the US dollar: 1, 5, 10, 25, 50 cent coins.
We choose a set of coins as
\[
S_{\mathrm{cent}}=\begin{pmatrix}
1 & 5 & 10 & 25 & 50\\
4 & 1 & 2 & 1 & 1
\end{pmatrix},
\]
so that each money amount from 1 to $\Amount{S_{\mathrm{cent}}}=104$
cents can be made by some coins of $S_{\mathrm{cent}}$.
The set of possible division is shown in Fig.~\ref{fig5},
which is obtained by enumerating all division directly.
A hierarchical structure can be seen but not perfectly.
Breaking of the hierarchical structure is due to the property
that some money amounts can be made in two ways,
e.g., $25=25\times1=5\times1+10\times2$
and $50=50\times1=5\times1+10\times2+25\times1$.

\section{Generalization of the number of players}
In this section, the number of players is changed.
For simplicity, we mainly deal with the binary coin set $S_{2,m}$
as in the most basic type of the problem stated in \secref{\ref{sec2}}.

We start from the coin dividing between four players $A$, $B$, $C$, and $D$.
Each way of division, written as a quadruplet $(n_A, n_B, n_C, n_D)$,
gives a point in a four-dimensional space,
and we represent the point as
$(n_A, n_B, n_C, n_D)
=n_A\Vector{e}_A+n_B\Vector{e}_B+n_C\Vector{e}_C+n_D\Vector{e}_D$.
Let us denote by $\Triangle_m^{(4)}$
the point set corresponding to the possible division of the coins.
The superscript ``$(4)$'' represents explicitly the number of players.

By using an argument similar to that in \secref{\ref{sec2}},
$\Triangle_m^{(4)}$ is determined inductively as
\begin{align*}
\Triangle_1^{(4)}&=\{\Vector{e}_A, \Vector{e}_B, \Vector{e}_C, \Vector{e}_D\},\\
\Triangle_m^{(4)}&=(\Triangle_{m-1}^{(4)}+2^{m-1}\Vector{e}_A)\cup
		(\Triangle_{m-1}^{(4)}+2^{m-1}\Vector{e}_B)\cup
		(\Triangle_{m-1}^{(4)}+2^{m-1}\Vector{e}_C)\cup
		(\Triangle_{m-1}^{(4)}+2^{m-1}\Vector{e}_D).
\end{align*}
The first equation signifies that
$\Triangle_1^{(4)}$ consists of four points
corresponding to who receives the coin of face value 1.
$\Triangle_1^{(4)}$ forms the vertices
of a regular tetrahedron whose sides are
$\|\Vector{e}_A-\Vector{e}_B\|=\|\Vector{e}_A-\Vector{e}_C\|
	=\cdots=\|\Vector{e}_C-\Vector{e}_D\|=\sqrt{2}$.
The second equation signifies that
$\Triangle_m^{(4)}$ is made up of four subsets
corresponding to who receives the coin of face value $2^{m-1}$,
which generates a nested structure of tetrahedra.

\begin{figure}[!t]\centering
\mbox{\raisebox{45mm}{(a)}}
\includegraphics[clip]{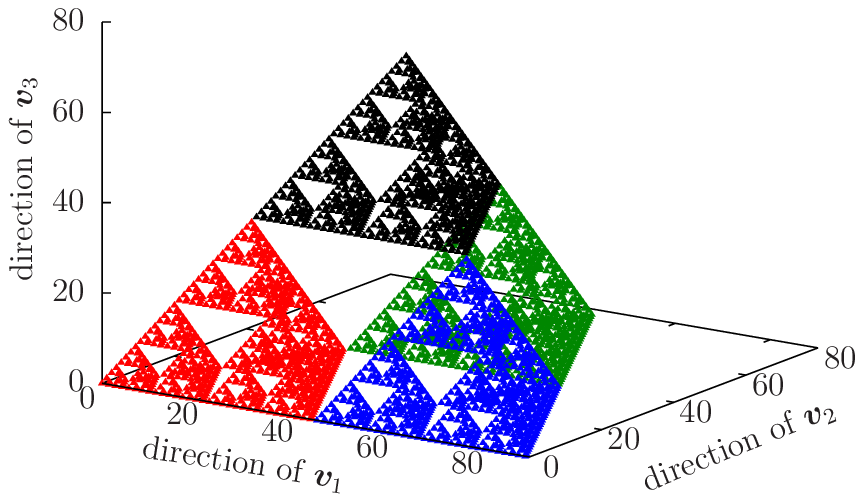}
\hspace{5mm}
\mbox{\raisebox{45mm}{(b)}}
\mbox{\raisebox{20mm}{\includegraphics[clip]{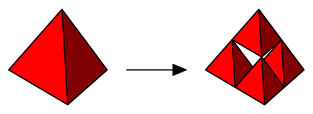}}}
\caption{
(a) A three-dimensional illustration of $\Triangle_m^{(4)}$ with $m=6$
using mapping \eqref{eq:GramSchmidtmap}.
It forms a Sierpinski tetrahedron.
(b) The iteration scheme of the Sierpinski tetrahedron.
}
\label{fig6}
\end{figure}

$\Triangle_m^{(4)}$ is a point set in the four-dimensional space,
but it lies on a three-dimensional affine hyperplane
given by $n_A+n_B+n_C+n_D=2^m-1$.
Hence we can visualize $\Triangle_m^{(4)}$
by taking a three-dimensional coordinate suitably.
Setting three vectors
$\Vector{u}_1:=\Vector{e}_A-\Vector{e}_D$,
$\Vector{u}_2:=\Vector{e}_B-\Vector{e}_D$, and
$\Vector{u}_3:=\Vector{e}_C-\Vector{e}_D$,
we get
\begin{align*}
n_A\Vector{e}_A+n_B\Vector{e}_B+n_C\Vector{e}_C+n_D\Vector{e}_D
&=n_A(\Vector{e}_A-\Vector{e}_D)+n_B(\Vector{e}_B-\Vector{e}_D)
	+n_C(\Vector{e}_C-\Vector{e}_D)+(2^m-1)\Vector{e}_D\\
&=n_A\Vector{u}_1+n_B\Vector{u}_2+n_C\Vector{u}_3+(2^m-1)\Vector{e}_D.
\end{align*}
That is,
$(n_A,n_B,n_C,n_D)\in\Span\{\Vector{u}_1,\Vector{u}_2,\Vector{u}_3\}+(2^m-1)\Vector{e}_D$.
One can find an orthonormal basis $\{\Vector{v}_1,\Vector{v}_2,\Vector{v}_3\}$
of subspace $\Span\{\Vector{u}_1, \Vector{u}_2, \Vector{u}_3\}$
by the Gram-Schmidt process, as
\[
\Vector{v}_1=\frac{\Vector{u}_1}{\sqrt{2}},\quad
\Vector{v}_2=\frac{2\Vector{u}_2-\Vector{u}_1}{\sqrt{6}},\quad
\Vector{v}_3=\frac{3\Vector{u}_3-\Vector{u}_2-\Vector{u}_1}{2\sqrt{3}}.
\]
Therefore, we have an isometric embedding
\begin{equation}
\Triangle_m^{(4)}\ni(n_A,n_B,n_C,n_D)\mapsto
\left(\frac{2n_A+n_B+n_C}{\sqrt{2}},\frac{3n_B+n_C}{\sqrt{6}},
	\frac{2n_C}{\sqrt{3}}\right)
\in\R^3.
\label{eq:GramSchmidtmap}
\end{equation}
The $i$-th component on the right-hand side is the scalar product
of $n_A\Vector{u}_1+n_B\Vector{u}_2+n_C\Vector{u}_3$ and $\Vector{v}_i$
($i=1,2,3$).
In Fig.~\ref{fig6} (a), we illustrate $\Triangle_m^{(4)}$ with $m=6$
embedded in a three-dimensional space
by using mapping \eqref{eq:GramSchmidtmap}.
$\Triangle_m^{(4)}$ is
a finite approximation of the three-dimensional Sierpinski gasket
or so-called Sierpinski tetrahedron, which is intuitively obtained
by iterating the removal process shown in Fig.~\ref{fig6} (b).
By analogy with the results in the previous section,
the change of a coin system affects
the number of layers of the Sierpinski tetrahedron.

The coin-dividing problem between $s(>4)$ players is considered as well.
Taking an orthonormal basis $\{\Vector{e}_1,\ldots,\Vector{e}_s\}$
of an $s$-dimensional space,
we have a recursive equation of $\Triangle_m^{(s)}$
\[
\Triangle_1^{(s)}=\{\Vector{e}_1,\ldots,\Vector{e}_s\},\quad
\Triangle_m^{(s)}
=\bigcup_{\sigma=1}^s(\Triangle_{m-1}^{(s)}+2^{m-1}\Vector{e}_\sigma).
\]
$\Triangle_m^{(s)}$ becomes
$(s-1)$-dimensional Sierpinski gasket as $m$ increases,
but we cannot visualize such a high-dimensional object.

Of course,
we can treat of the coin-dividing problem between one or two players.
For the dividing by one player $A$, there is only one {\it trivial} dividing
of which all coins go to $A$;
therefore, $\Triangle_m^{(1)}$ is a one-point set.
As for the dividing between two players $A$ and $B$,
the solution is
$\Triangle_m^{(2)}=\{(n_A,n_B)|n_A=0,1,\ldots,2^m-1, n_A+n_B=2^{m}-1\}$,
consisting of $2^m$ points which are equally spaced.
The inductive relation in this case is
\[
\Triangle_m^{(2)}=(\Triangle_{m-1}^{(2)}+2^{m-1}\Vector{e}_A)\cup
	(\Triangle_{m-1}^{(2)}+2^{m-1}\Vector{e}_B),
\]
which is converted to
\[
2^{-m}\Triangle_m^{(2)}=f_A(2^{-(m-1)}\Triangle_{m-1})\cup f_B(2^{-(m-1)}\Triangle_{m-1}),
\]
where $f_A$ and $f_B$ are given in Eq.~\eqref{eq:fP}.
It is easy to find that
the attractor of IFS $\{f_A,f_B\}$ is a line segment
$\{(n_A,n_B)| n_A\ge0, n_B\ge0, n_A+n_B=1\}$.
In conclusion, nontrivial fractal pattern is not observed in these two cases.
However, we explain in the next section
that coin-dividing between two players with a suitable coin set
can generate the Cantor set.

\section{Cantor set between two players}
In the previous section,
no interesting fractal structure emerges in coin dividing between two players.
We see in this section that the Cantor set is obtained
if the coin set is not complete.

Let us study on trial the coin-dividing problem
by two players $A$ and $B$ with the set of coins
\[
S_{3,1,m} =
\begin{pmatrix}1 & 3 & 9 &\cdots& 3^{m-1}\\
	1 & 1 & 1&\cdots& 1\end{pmatrix}.
\]
Some money amounts,
e.g., 2, 5, 6, 7, and 8, cannot be made using these coins---
this {\it incompleteness} is a different property from $S_{2,m}(=S_{2,1,m})$.
We let $\Cantor_m$ denote the set of pairs $(n_A, n_B)$ of dividing coins $S_{3,1,m}$.

\begin{figure}[!t]\centering
\mbox{\raisebox{47mm}{(a)}}
\includegraphics[clip]{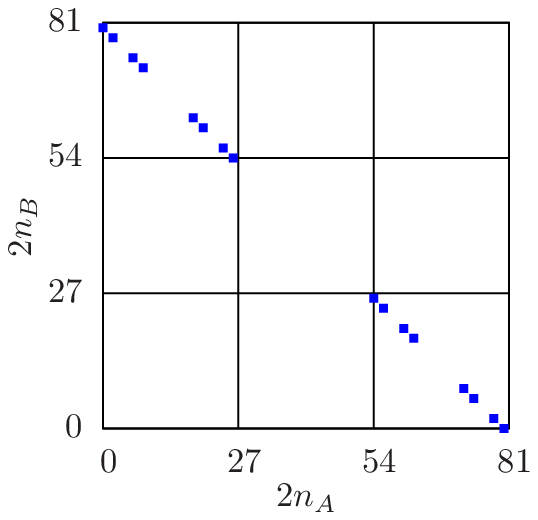}
\hspace{5mm}
\mbox{\raisebox{47mm}{(b)}}
\includegraphics[clip]{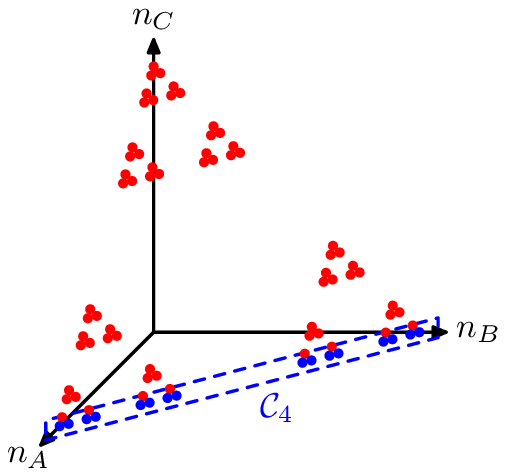}
\caption{
(a) The possible division of $S_{3,1,4}$.
For ease of display, two axes are chosen to $(2n_A,2n_B)$.
(b) The intersecting points of
$\Triangle_{3,1,4}$ and the $(n_A,n_B)$-plane is $\Cantor_4$.
}
\label{fig7}
\end{figure}

The inductive relation of $\Cantor_m$ is
\[
\Cantor_m=(\Cantor_{m-1}+3^{m-1}\Vector{e}_A)
	\cup(\Cantor_{m-1}+3^{m-1}\Vector{e}_B),
\]
and is equivalently
\[
3^{-m}\Cantor_m=\left(\frac{3^{-(m-1)}\Cantor_{m-1}+\Vector{e}_A}{3}\right)
\cup\left(\frac{3^{-(m-1)}\Cantor_{m-1}+\Vector{e}_B}{3}\right)
=g_A(3^{-(m-1)}\Cantor_{m-1})\cup g_B(3^{-(m-1)}\Cantor_{m-1}),
\]
where $g_A(\Vector{x})=(\Vector{x}+\Vector{e}_A)/3$ and
$g_B(\Vector{x})=(\Vector{x}+\Vector{e}_B)/3$ are contraction maps.
The set $3^{-m}\Cantor_m$ converges to
the attractor of IFS $\{g_A, g_B\}$,
and as stated in \secref{\ref{sec1}}, it is the Cantor set
spanned between $\Vector{e}_A$ and $\Vector{e}_B$.
In addition,
$\Cantor_m$ can be regarded as a subset of $\Triangle_{3,1,m}$
formed by the points of $n_C=0$,
or the intersecting points of $\Triangle_{3,1,m}$ and the $(n_A,n_B)$ plane
(see Fig.~\ref{fig7} for reference).

Let us study more about the coin-dividing problem and the Cantor set.
We introduce a mapping
$\varphi_m:\Cantor_m\ni(n_A,n_B)\mapsto 2\cdot3^{-m}n_A\in [0,1]$.
$\varphi_m$ is the composite mapping of
a projection $\pi:(x,y)\mapsto x$ and
a scaling function $\rho_m:x\mapsto 2x/3^m$.
By definition, money amount $n_A$
is written as $n_A=\sum_{k=0}^{m-1}\chi_k 3^k$,
where $\chi_k\in\{0,1\}$ is the indicator of
whether the coin of face value $3^k$ goes to $A$ or not.
Thus,
\[
\varphi_m(\Cantor_m)
=\left\{\left.\sum_{k=0}^{m-1}\frac{2\chi_k}{3^{m-k}}
	\right|\chi_k\in\{0,1\}\right\}.
\]
$\varphi_m(\Cantor_m)$ consists of
numbers within $[0,1]$ whose base-3 representations are up to $m$ digits
with entirely 0s and 2s,
and formally $\lim_{m\to\infty}\varphi(\Cantor_m)\subset[0,1]$ consists of
numbers whose base-3 representation are entirely 0s and 2s.
As is known well \cite{Peitgen},
this is identical with another definition of the Cantor set.

Generalization to the coin set
\[
S_{r,1,m}=\begin{pmatrix}1 & r & r^2 &\cdots& r^{m-1}\\
	1 & 1 & 1&\cdots& 1\end{pmatrix}
\]
is straightforwardly.
By way of the inductive relation for the set $\Cantor_{r,m}$ of division,
$r^{-m}\Cantor_{r,m}$ converges to a Cantor-like fractal
whose Hausdorff dimension is
\[
D_r=\frac{\ln2}{\ln r}.
\]
In particular, $D_3=\ln2/\ln3\approx0.631$ is the dimension of the Cantor set.
Moreover,
$D_2=\ln2/\ln2=1$ is consistent with
the result that the attractor of coin-dividing of $S_{2,m}(=S_{2,1,m})$
becomes a line segment, as stated at the end of the previous section.

\section{Conclusion}
In this article,
we have treated of the {\it coin-dividing problem}
which seeks the set of the possible division of a set of coins.
By considering an appropriate scaling limit,
the set of points in dividing of coins $S_{r,c,m}$
converges to a fractal set as $m$ tends to infinity.
This result follows from
a remarkable property that
increment of the coin types $m$ represents
the action of an iterated function system.
The parameters $r$ and $c$ of $S_{r,c,m}$ are related to
the contraction ratio and the number of maps of the iterated function system,
respectively.
Depending on the magnitude of $c$ and $r-1$,
the coin-dividing fractal belongs to one of three classes:
totally disconnected if $c<r-1$,
finitely ramified if $c=r-1$,
and infinitely ramified if $c>r-1$.
In particular,
we have obtained the Sierpinski gasket when three players divide $S_{2,1,m}$,
and the Cantor set when two players divide $S_{3,1,m}$.

\section*{Acknowledgments}
The author is very grateful to Dr.~Makoto Katori
for fruitful and instructive discussion.
The author is also grateful to Dr.~Yoshihiro Yamazaki and
Dr.~Jun-ichi Wakita for their comments and discussion.

\end{document}